\documentclass{osa-article}

\journal{osajournal}

\articletype{Research Article}

\usepackage{tabularx,amsmath,amsfonts,boxedminipage,graphicx,mathtools,appendix,physics}
 \usepackage{epstopdf, float}
 \usepackage{algorithm2e}
 \usepackage{cite}
 \usepackage{subfigure}
 \newenvironment{remark}[1][Remark]{\begin{trivlist}
 \item[\hskip \labelsep {\bfseries #1}]}{\end{trivlist}}
 \numberwithin{equation}{section}
 \numberwithin{figure}{section}
 \RestyleAlgo{boxruled}

\begin{document}
\title{An Optimal Control Approach to Gradient-Index Design for Beam Reshaping}

\author{J. Adriazola,\authormark{1,*} and R. H. Goodman,\authormark{1}}

\address{\authormark{1}Department of Mathematical Sciences and Statistics, New Jersey Institute of Technology, University Heights, Newark, New Jersey, 07102, USA}

\email{\authormark{*}ja374@njit.edu} 

\begin{abstract}
We address the problem of reshaping light in the Schrödinger optics regime from the perspective of optimal control theory. In technological applications, Schrödinger optics is often used to model a slowly-varying amplitude of a para-axially propagating electric field where the square of the waveguide's index of refraction is treated as the potential. The objective of the optimal control problem is to find the controlling potential which, together with the constraining Schrödinger dynamics, optimally reshape the intensity distribution of Schrödinger eigenfunctions from one end of the waveguide to the other. This work considers reshaping problems found in work due to Kunkel and Leger, and addresses computational needs by adopting tools from the quantum control literature. The success of the optimal control approach is demonstrated  numerically.
\end{abstract}

\section{Introduction}
Humans have been reshaping light for thousands of years, and it remains an active research area to this day,
from the ancient Assyrians' introduction of primitive lenses circa 750 B.C.E~\cite{Nimrud} to designs based on the sophisticated techniques of optimal transport~\cite{Froese}.
Requiring a laser beam to have a specified irradiance distribution has diverse and broad applications which include laser/material processing, laser/material interaction studies, fiber injection systems, optical data image processing, and lithography~\cite{Dickey}. Geometric optics is the simplest physical setting in which to study beam reshaping, and one that is often chosen. However, in the presence of diffractive effects, the wave nature of light must be accounted for, as is often the case in nano-scale optical technologies. 

\begin{figure}[htbp]
\begin{centering}
{\includegraphics[width=1\textwidth]{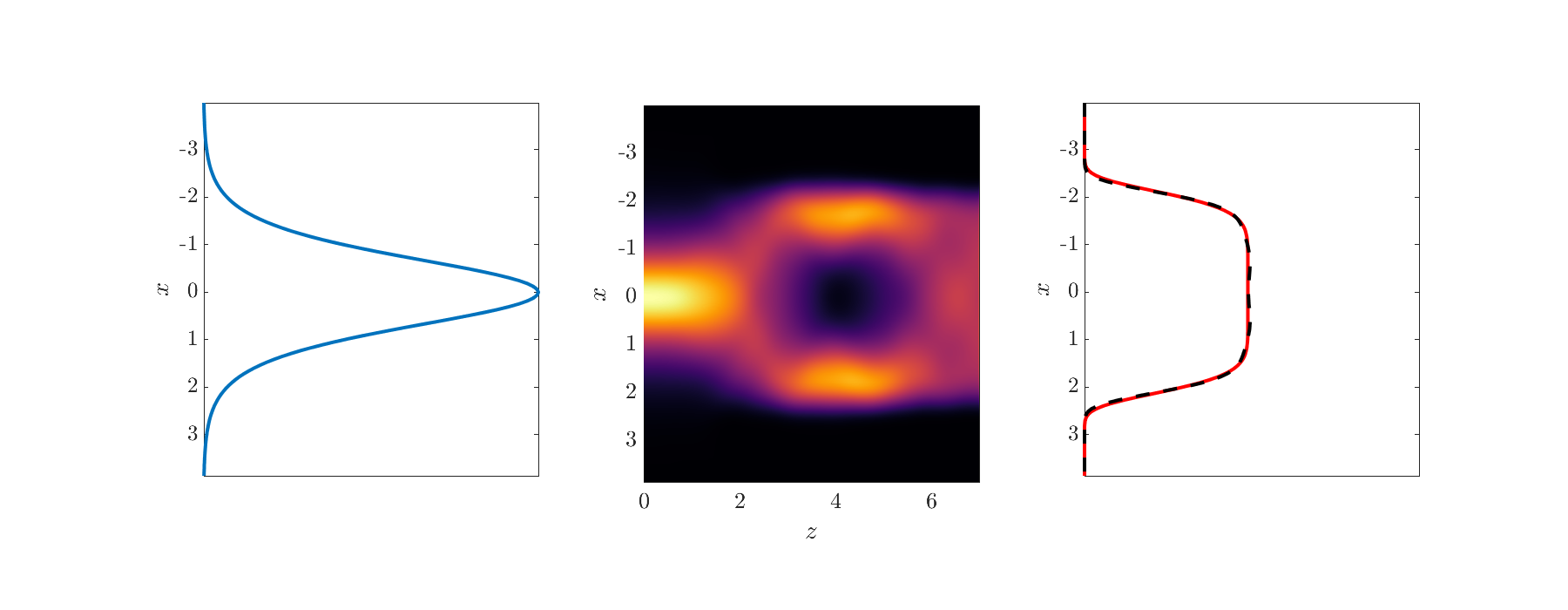}}
\caption{An example of reshaping light with a peaked intensity profile into one with a more uniform profile. This paper details the methods used to find such a mapping.}
\label{fig:introsquare}
\end{centering}
\end{figure}

\begin{figure}[htbp]
\begin{centering}
{\includegraphics[width=1\textwidth]{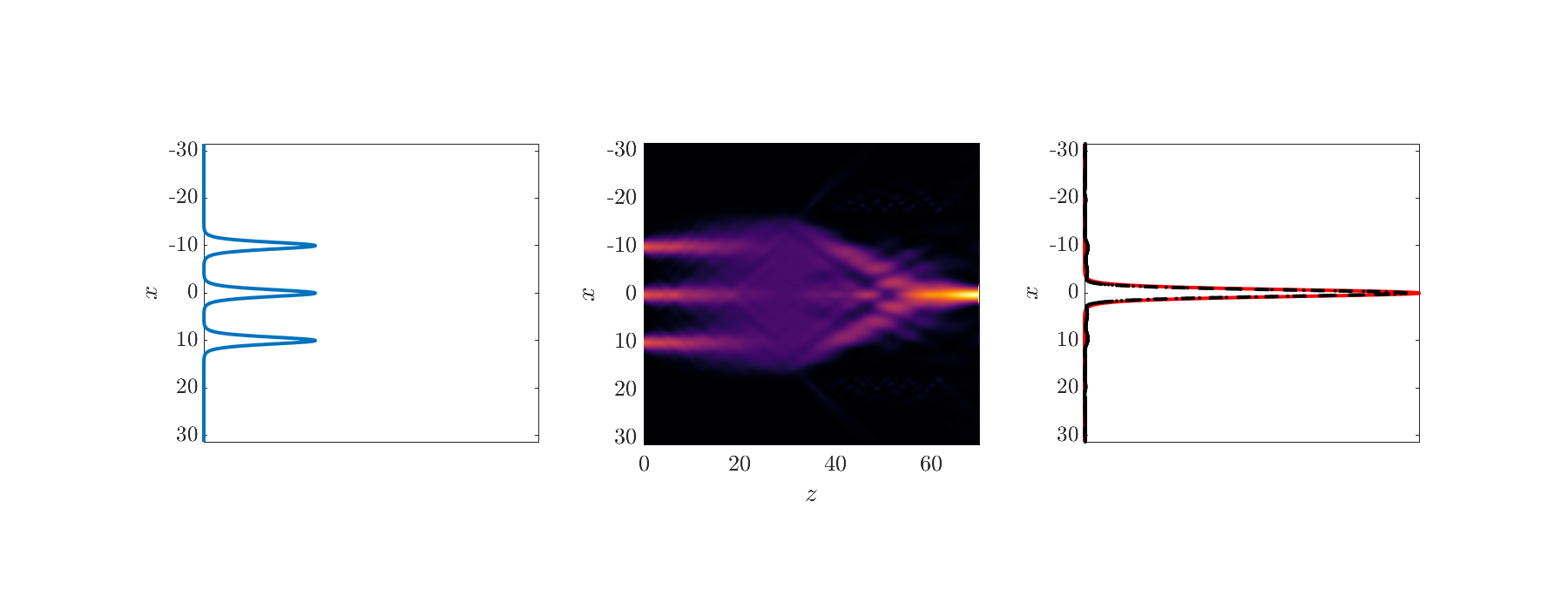}}
\caption{An example of light reshaping in which three pulses of light are combined into one using optimal control theory. More detail about the construction of the GRIN component which combines these pulses is provided throughout the paper.}
\label{fig:introcombine}
\end{centering}
\end{figure}

In recent work, Kunkel and Leger~\cite{Kunkel,MinsterKunkel:20} successfully reshape laser beams in the presence of diffraction. They demonstrate that the phase retrieval method~\cite{Fienup} is a viable means for numerically constructing a gradient-index (GRIN) optical waveguide which reshapes light into an intended intensity distribution. Figure~\ref{fig:introsquare} shows our computation of an example from~\cite{Kunkel} in which light is transformed from a sharply peaked intensity profile to nearly uniform one. Figure~\ref{fig:introcombine} shows the application of our methods to another example from~\cite{Kunkel} in which the GRIN combines multiple localized intensity distributions into one.

A disadvantage of the phase retrieval method is that generalizing it to either higher spatial dimensions or generalizing the dynamical constraints may be difficult. Indeed, despite achieving great success, Kunkel and Leger show several necessary adjustments must be made in order to adapt their previous methodology in two spatial dimensions~\cite{Kunkel} to the case of three~\cite{MinsterKunkel:20}. On the other hand, optimal control theory, an extension of the calculus of variations~\cite{Gelfand,Calculus1989}, provides a more general alternative method. The chief advantage of using optimal control theory is in its abstract framework which easily handles entire classes of optimization problems at once, independent of its dimension or class of constraints.

In this work, we pose an optimal control problem with an objective functional, first used in the context of high-fidelity quantum fluid manipulations by Hohenester, et al.~\cite{Hohenester}, constrained by the following standard model for paraxial light beam propagation. Consider an electromagnetic field propagating transversely through a linear waveguide, i.e., a waveguide through which the electrical field responds linearly to the polarization of the propagation media. Asssume the propagating field is time-harmonic, has negligible magnetic field components, and satisfies the hypothesis of the paraxial approximation, namely, that the direction of propagation does not deviate signficiantly from the axial direction defined by the waveguide. Then, one can show that Schrödinger's equation, in dimensionless form,
\begin{equation}\label{eq:Intro2}
    i\psi_z=-\frac{1}{2}\triangle\psi+V(x,z)\psi,
\end{equation}
arises as a slowly varying amplitude approximation to the variable-coefficient Helmholtz equation~\cite{JWGood}.

Here, $z$ is the axis of propagation, $x$ is the transverse direction, $\triangle$ is the Laplacian in the transverse direction, $V(x,z)$ is proportional to the square of a spatially varying refractive index, and the wavefunction $\psi(x,z)$ is interpreted as a spatially varying complex electric amplitude.  We assume the propagation media is lossless, hence the potential $V$ is a real function of the waveguide coordinates. The paraxial approximation is often studied because the numerical solution of Schrödinger's equation is significantly cheaper computationally, and easier to understand analytically, than either the full numerical solution of  Helmholtz's or Maxwell's equations. 

In posing the design problem, we make the simplifying assumption that the potential can be written in the form $V(x,z)=V(x,u(z))$ where $V(x,u)$ is a one-parameter family of potentials. Thus the design of a reshaping potential $V(x,z)$ is reduced to a search for a one-dimensional optimal control $u(z)$. The light reshaping problem in this paper is therefore: Find the optimal control $u(z)$ that best transforms the intensity distribution of an initial Schrödinger state $\varphi_0(x)$
into the intensity distribution of the desired state $\varphi_d(x)$ satisfying
\begin{subequations}\label{eq:initdes}
\begin{align}
-\frac{1}{2}\triangle\varphi_0(x)+V(x,u(0))\varphi_0(x)&=\lambda_0\varphi_0(x),\\
-\frac{1}{2}\triangle\varphi_d(x)+V(x,u(l))\varphi_d(x)&=\lambda_l\varphi_d(x),
\end{align}
\end{subequations}
i.e., the initial and desired states are eigenfunctions, of the time-independent Schrödinger operator $P=-\frac{1}{2}\triangle+V(x,u(z))$, at $z=0$ and at the end of a specified propagation length $l$, respectively. Thus, we formulate the problem of designing an optimal coupler between two waveguides with different transverse profiles and their eigenpairs $\left(\varphi_0,\lambda_0\right)$ and $\left(\varphi_d,\lambda_l\right)$.

\subsection{Structure of the Paper}
In Section~\ref{section:beamprob}, we precisely state the eigenfunction reshaping problem considered throughout this work. The problem is similar to quantum optimal control problems previously considered  in the literature, e.g.~\cite{Hohenester, Mennemann, Nature}. We discuss, in detail, our assumptions about the control problem, and provide the optimality conditions given by the Euler-Lagrange equations~\cite{Gelfand}. 


In Section~\ref{section:beamnum}, we provide an overview of the numerical methods used in solving the control problem posed in Section~\ref{section:beamprob}. The procedure is a combination of a global, non-convex method followed by a local, iterative method. In the context of numerical optimal control, this approach is called a $\mathbf{hybrid\ method}$~\cite{Sorensen}. Hybrid optimization methods, when used appropriately, can overcome non-convexity, yet still remain computationally efficient. 

The group of Calarco was the first to use the type of non-convex method we use in this work~\cite{Doria,Caneva}. This method reduces the dimensionality of the control problem so that standard global search routines based on stochastic optimization can be used. Since stochastic methods come at the cost of slow convergence near local minima~\cite{BoydV}, local methods are then used to accelerate convergence toward the nearest minimum. The local method we use is a gradient descent, due to Von Winckel and Borzi~\cite{Borzi}, called GRAPE, which ensures that controls remain in the admissible search spaces used throughout this work.

In Section~\ref{section:problems}, we address many of the practical and computational aspects arising from the specific beam reshaping problems of interest. We use reductions which greatly simplify the computational complexity of the problem, and greatly aid in efficiently searching the space of reshaping potentials. The success of these reductions, together with the methods detailed in Section~\ref{section:beamnum}, is demonstrated numerically on the two reshaping problems previously shown in Figures~\ref{fig:introsquare} and~\ref{fig:introcombine}.

\subsection{Notation and Conventions}
We make use of various function spaces when stating the optimal control problem. For example, general Banach spaces are denoted by $\mathcal{B}$. The Lebesgue space denoted by $L^p(\Omega)$, where $\Omega$ is a measurable set, is the equivalence class of measurable functions which agree almost everywhere such that the norm
\begin{equation}
    \|f\|_{L^p(\Omega)}:=\left(\int_{\Omega}|f|^pd\mu\right)^{\frac{1}{p}}
\end{equation}
is finite. Similarly, the Sobolev space $H^{k}(\Omega)$ is the space of $k$-times weakly differentiable functions $f$, with respect to $x\in\Omega$, whose norm,
\begin{equation}\label{eq:Sob}
    \|f\|_{H^{k}(\Omega)}:=\left(\sum_{j=0}^{k}\left|\left|\partial^j_xf\right|\right|^2_{L^2(\Omega)}\right)^{\frac{1}{2}}
\end{equation}
is finite. The space of essentially bounded functions $L^{\infty}(\Omega)$ is the space where
\begin{equation}
    \|f\|_{L^{\infty}(\Omega)}:=\mathrm{ess}\sup_{x\in\Omega}|f(x)|<\infty.
\end{equation}

Homogeneous Sobolev spaces, denoted by $\dot{H}^{k}(\Omega)$, are the spaces of functions such that $\left\|\partial_x^kf\right\|_{L^{2}(\Omega)}$ is finite. A traceless Sobolev space, denoted ${H}_0^{k}(\Omega)$, is the space of functions in ${H}^{k}(\Omega)$ which vanish on the boundary $\partial\Omega$. The space of $k$-times continuously differentiable functions is denoted $C^k(\Omega)$, and the space of essentially bounded $C^k(\Omega)$ functions is denoted by
\begin{equation}
    C^k_b(\Omega):=C^k(\Omega)\cap L^{\infty}(\Omega).
\end{equation}

The notation $\mathcal{B}_1(\Omega_1;\mathcal{B}_2(\Omega_2))$ is understood as the space of functions $f$ such that $\|f(\Omega_1,\cdot)\|_{\mathcal{B}_2(\Omega_2)}\in{\mathcal{B}_1(\Omega_1)}$. Spaces where each element is compactly supported on $\Omega$ are denoted by $\mathcal{B}_c(\Omega)$. Lastly, the notation $^{\dag}$ denotes Hermitian conjugation.

\section{Optimal Control Framework}\label{section:beamprob}
The salient elements of the problem structure we consider are due to Hohenester, et al.~\cite{Hohenester}, which uses the following objective functional
\begin{equation}\label{eq:BeamObj}
J=\frac{1}{2}\left(\|\varphi_d(\cdot)\|_{L^2(\mathbb{R}^n)}^4-\left|\left\langle \varphi_d(\cdot),\psi(\cdot,l)\right\rangle\right|^2_{L^2(\mathbb{R}^n)}\right)+\frac{\gamma}{2}\int_0^l\left|\partial_z u\right|^2dt,
\end{equation}
where $\gamma>0$ and $z\in(0,l)$ is the axial coordinate, with $l>0$.

The objective functional $J$ involves the \emph{infidelity} 
\begin{equation}\label{eq:Jinfed}
    J_{\rm infidelity}= \frac{1}{2}\left(\|\varphi_d(\cdot)\|_{L^2(\mathbb{R}^n)}^4-\left|\left\langle \varphi_d(\cdot),\psi(\cdot,l)\right\rangle\right|^2_{L^2(\mathbb{R}^n)}\right)
\end{equation} 
which penalizes misalignments of the computed function $\psi(x,l)$ with respect to the desired state $\varphi_d(x)$. In the language of optimal control theory~\cite{Bryson,Calculus1989}, the infidelity is called a terminal cost. This objective functional disregards the physically unimportant global phase difference between the desired and computed states,  a significant advantage over the typical least-squares approach.

The second contribution to the objective, the running cost over $[0,l]$, is a \emph{regularization} of the control function $u(z)$. This penalizes the use of control functions with large ${\dot{H}}([0,l])$ norms, and is well-known in the literature as a type of Tikhonov regularization~\cite{Tikhonov1995}. The introduction of this regularization conditions the optimal control problem. Indeed, Hintermuller, et al., prove the control framework of Hohenester, et al., is well-posed with the introduction of a Tikhonov regularization, i.e. there exists a control $u\in H^1([0,l])$ that minimizes the objective $J$~\cite{Hintermuller}.

The optimal control problem we consider in this paper is the following:
\begin{equation}\label{eq:BeamProb}
    \inf_{u\in\mathcal{U}}J
\end{equation}
subject to Schrödinger's equation~\eqref{eq:Intro2} with the initial and desired states $\varphi_0$ and $\varphi_d$ satisfying Equations~\eqref{eq:initdes}.
The search for optimal controls is performed over the admissible class $\mathcal {U}=\left\{u\in H^1\left([0,l]\right):u(0)=u_0,u(l)=u_l\right\}$. We assume the eigenfunctions $\varphi_d$ and $\varphi_0$ are both in the space $H^{1}(\mathbb{R}^n)$. We also assume that the eigenfunctions $\varphi_0$ and $\varphi_d$ have unit intensity, i.e., $\|\varphi_0\|_{L^2({\mathbb{R}^n})}=\|\varphi_d\|_{L^2({\mathbb{R}^n})}=1$, so that the infimum of the infidelity~\eqref{eq:Jinfed} is 0. Lastly, we assume the potential $V(x,u(z))$ is in the space $C_b^{0}([0,l];H^{1}\left(\mathbb{R}^n\right))$ for every $u\in\mathcal{U}$.

\begin{remark}
With the above assumptions in place, the regularity of the wavefunction $\psi$ solving Equation~\eqref{eq:Intro2} is known~\cite{MASPER}; $\psi\in C^1([0,l];H^1(\mathbb{R}^n))$. Moreover, the control problem with objective functional~\eqref{eq:BeamObj} is well-posed for sufficiently large $\gamma>0$~\cite{Hintermuller}. 
\end{remark}

By letting
\begin{equation*}
    J\to J+\int_0^l\int_{\mathbb{R}^n} p^{\dag}\left(i\psi_z+\frac{1}{2}\triangle\psi-V(x,z)\psi\right)dz
\end{equation*}
where $p$ is a Lagrange multiplier, and using standard arguments from the calculus of variations~\cite{Gelfand,Witelski}, it is straightforward to show that the optimality conditions of Problem~\eqref{eq:BeamProb} are given by
\begin{subequations}\label{eq:optcond}
\begin{alignat}{2}
\label{eq:stateeq}
i\partial_z\psi&=-\frac{1}{2}\triangle\psi+V(x,z)\psi,\qquad &\psi(x,0)=\varphi_0(x), \\
\label{eq:costateeq}
i\partial_zp&=-\frac{1}{2}\triangle p+V(x,z) p,\qquad &ip(x,l)=\left\langle \varphi_d,\psi(x,l)\right\rangle_{L^2(\mathbb{R}^n)}\varphi_d, \\
\label{eq:controleq}
\gamma\partial_z^2{u}&=-\Re\left\langle p,\partial_uV\psi\right\rangle_{L^2(\mathbb{R}^n)},\qquad &u(0)=u_0,\ u(l)=u_l.
\end{alignat}
\end{subequations}
Equation~\eqref{eq:costateeq} is the adjoint equation of Equation~\eqref{eq:stateeq} and governs the axial evolution of the Lagrange, or costate, multiplier $p$ backwards from its terminal condition at $z=l$. The similarity of Equation~\eqref{eq:stateeq} and Equation~\eqref{eq:costateeq} is due to the self-adjoint nature of the Schrödinger operator $P= -\frac{1}{2}\triangle +V(x,z)$.  Equation~\eqref{eq:controleq} governs the optimal control $u$, and together with the boundary conditions defined through the admissible class $\mathcal{U}$, is a boundary value problem on $[0,l]$. 

Equations~\eqref{eq:stateeq} and~\eqref{eq:stateeq} are both solved via a second-order Fourier split-step method, where the $z-$dependence of the potential is handled by the midpoint method. We also note that Equation~\eqref{eq:controleq} will not be solved numerically, but will instead be reinterpreted in the context of the optimization method discussed in Section~\ref{section:Local}.

Consider the so-called reduced objective functional
\begin{equation}\label{eq:reducedcost}
\mathcal{J}:\mathcal{U}\to\mathbb{R},\quad u\mapsto\mathcal{J}[u]:=J\left[\psi(u),u\right].
\end{equation}
Let $u^*$ denote an optimal control, and define $\psi^*:=\psi(u^*)$, $p^*:=p(u)$. Since the optimal control problem~\eqref{eq:BeamProb} is well-posed, then for every $u\in\mathcal{U}$,
\begin{equation}\label{eq:optineq}
    \mathcal{J}[u]\geq \mathcal{J}[u^*]=\min_{u\in\mathcal{U}}\mathcal{J}
\end{equation}
if and only if the optimal triple $(\psi^*,p^*,u^*)$ satisfies
Equations~\eqref{eq:optcond}. For this reason, pursuing numerical approximations of Equations~\eqref{eq:optcond} and the optimality condition~\eqref{eq:optineq} when searching for the optimal control $u^*$ is meaningful.

\section{Numerical Optimization Methods}\label{section:beamnum}

In order to solve Problem~\eqref{eq:BeamProb}, we use a \emph{hybrid} optimization method; a combination of a global, non-convex method followed by a local, iterative method. The methodology we use in this paper is similar to one used by S{\o}rensen, et al.~\cite{Sorensen}, and allows for the use of a global search routine based on stochastic optimization to overcome non-convexity. Non-convex objective functions may, of course, possess many local minima, and a global method seeks to efficiently search for a near-optimal one. By then feeding results from the global method into the local one, convergence near the local minimum is accelerated. We previously used this methodology in~\cite{me}, and more specific details about the numerical optimization is provided there.

\subsection{The Global Method}\label{section:Global}
The first step in the hybrid method is to use a Galerkin method which reduces the complexity of the optimal control problem so that standard non-convex nonlinear programming (NLP) techniques can be applied. This step relies on choosing controls from the span of an appropriately chosen finite set of basis functions so that the optimization is performed over a relatively small set of unknown coefficients. The choice of basis is such that controls remain in the appropriate admissible space $\mathcal{U}$ in the context of the control problem~\eqref{eq:BeamProb}. 

We choose to use the representation
\begin{equation}\label{eq:randclass}
    u_r(z)=\mathcal{P}(z;u_0,u_l,l)+\sum_{j=0}^{N-1}\varepsilon_j\varphi_j(z;l),\quad z\in[0,l],\\
\end{equation}
where $\mathcal{P}$ is a fixed function satisfying the boundary conditions of Equation~\eqref{eq:controleq}, $\varphi_j(z)$ is a basis function with vanishing boundary conditions, and the coefficients $\varepsilon_j$ are parameters to be optimized over. It is clear that if the polynomial $\mathcal{P}$ and the basis functions $\varphi_j$ are chosen well enough, then control ansatz~\eqref{eq:randclass} reliably simplifies the optimal control problem. An effective Galerkin approximation must be constructed the set of basis functions $N$ simultaneously large enough to define an accurate approximation, yet small enough so that the overall procedure remains computationally inexpensive. In this work, we use 15 basis functions. This reduces the  optimization problem to a small-scale NLP problem that can be solved using standard techniques. 

To solve the resulting NLP problem, we use differential evolution (DE)~\cite{Storn}. DE is a stochastic optimization method used to search for candidate solutions to non-convex optimization problems. The idea behind DE is a so-called genetic algorithm that draws inspiration from evolutionary genetics. DE searches the space of candidate solutions by initializing a population set of vectors, known as agents, within some chosen region of the search space. These vectors are then randomly mutated into a new population set, or generation. The mutation operates via two mechanisms: a weighted combination and a "crossover" which randomly exchanges "traits", or elements, between agents. 

DE ensures that the objective functional $\mathcal{J}$ decreases monotonically with each generation. As each iteration "evolves" into the next, inferior vectors "inherit" optimal traits from superior vectors via mutations. DE only allows mutations which are more optimal with respect to $\mathcal{J}$ to pass to the next generation. After a sufficient number of iterations, the best vector in the final generation is chosen as the candidate solution most likely to be globally optimal with respect to an objective functional.

Genetic algorithms, which require very few assumptions about the objective functional, are part of a wider class of optimization methods called metaheuristics. Although metaheuristics are useful for non-convex optimization problems, they do not guarantee about the global optimality of candidate solutions. Since the algorithm is stopped after a finite number of iterations, different random realizations return different candidate optimizers. For this reason, we use DE to search for candidate solutions and use these candidates in order to generate initial controls, through the representation~\eqref{eq:randclass}, for a method which guarantees local optimality up to some threshold.

\subsection{The Local Method}\label{section:Local}
We use a line search strategy due to Borzi and von Winckel called GRAPE~\cite{vonWinckel}. The GRAPE method is an appropriate generalization of the well-known gradient descent method from $\mathbb{R}^n$ to an appropriate affine function space which automatically preserves the boundary conditions of the admissible class $\mathcal{U}$ mentioned in the context of optimal control problem~\eqref{eq:BeamProb}. This method has been frequently applied in the quantum control literature; see for example~\cite{Hohenester,Mennemann,Sorensen}. 

To describe the GRAPE method, note that the optimal control problem~\eqref{eq:BeamProb} may be rewritten in the unconstrained form
\begin{equation}\label{eq:Lastbeamprob}
\min_{u\in \mathcal{U}}\mathcal{J}=
\min_{u\in \mathcal{U}}\int_0^l\mathcal{L}(\psi,\partial_z\psi,\partial_x^2{\psi},\psi^{\dag},p^{\dag},u,\partial_zu) dz,
\end{equation}
through routine manipulations of the objective~\eqref{eq:BeamObj} and use of the Lagrange multiplier $p(x,z)$.
The method of gradient descent, in this context, is given by following iteration
\begin{equation}\label{eq:graddesc}
u_{k+1}=u_k-\alpha_k\nabla_{u}\mathcal{L}
\big|_{u=u_k},
\end{equation} 
where the linear operator $\nabla_u$ is the gradient, or Fréchet derivative, of the Lagrangian $\mathcal{L}$ with respect to the control $u$. The stepsize $\alpha_k$ is chosen adaptively via the Armijo-Goldstein condition~\cite{BoydV}.  

Recall that the definition of a Fréchet derivative depends on the choice of function space in which it is to be understood. If the Fréchet derivative is understood in the sense of $L^2([0,l])$, then it can be identified with the functional derivative of the objective $\mathcal{J}$, which in this case can be shown to be
\begin{equation}
 \delta_u\mathcal{J}=-\gamma\partial_z^2{u}-\Re\left\langle p,\partial_uV\psi\right\rangle_{L^2(\mathbb{R}^n)}.
\end{equation}
This coincides with the Euler-Lagrange equation $\delta_u\mathcal{J}=0$ given by Equation~\eqref{eq:controleq}. If this choice is made, however, the increment $\alpha_k\nabla_{u}\mathcal{L}
\big|_{u=u_k}$ would not in general satisfy the boundary conditions on the control $u_k$, and the updated control $u_{k+1}$ would leave the admissible set $\mathcal{U}$. This problem is avoided by using a different function space $X$ defining the operator $\nabla_u$.

To this end, consider an arbitrary displacement $v\in C_c^{\infty}([0,l])$ and an arbitrary $\varepsilon>0$. We know Taylor's theorem holds, i.e., the series
\begin{equation}
J[u+\varepsilon v]=J[u]+\varepsilon\left\langle \nabla_{u}\mathcal{L}(u),v\right\rangle_{X}+\mathcal{O}(\varepsilon^2)
\end{equation}
holds term-by-term independently of the choice of the Hilbert space $X$ for sufficiently regular functionals $\mathcal{J}$. The GRAPE method chooses the function space $\dot{H}_0^1([0,l])$ for $X$.  By equating the directional, or Gateaux, derivatives with respect to $L^2([0,l])$ and with respect to $\dot{H}_0^1([0,l])$, we see that
\begin{align}\label{eq:weakform}
\begin{split}
&\left\langle \nabla_u \mathcal{L},v\right\rangle_{L^2([0,l])}=\left\langle \delta_u \mathcal{J},v\right\rangle_{L^2([0,l])}
\\
=&\left\langle \nabla_{u}\mathcal{L},v\right\rangle_{\dot{H}_0^1([0,l])}:=\int_0^l \partial_z\nabla_{u}\mathcal{L}\partial_zvdz
=-\left\langle \partial_z^2\nabla_{u}\mathcal{L},v \right\rangle_{L^2([0,l])},
\end{split}
\end{align}
where an integration by parts is used once along with the boundary conditions on $v$. 

Since this holds for all displacements $v\in C_c^{\infty}([0,l])$, we conclude, by the fundamental lemma of the calculus of variations~\cite{Gelfand}, the strong form of Equation~\eqref{eq:weakform}
\begin{equation} \label{eq:Poisson}
-\partial_z^2\nabla_{u}\mathcal{L} =\delta_u J,\  \ \nabla_u \mathcal{L}(0)=\nabla_u \mathcal{L}(l)=0,
\end{equation}
also holds. This renders an admissible gradient whose homogeneous Dirichlet conditions are induced by choosing increments specifically from the traceless space $\dot{H}_0^1([0,l])$. In order to solve the  boundary value problem~\eqref{eq:Poisson} for the control gradient $\nabla_u \mathcal{L}$, we use Chebyshev collocation~\cite{Trefethen}. 

\section{Beam Reshaping Problems}\label{section:problems}
\subsection{The Top Hat Problem}\label{section:THprob}
\begin{figure}[htbp]
    \centering
    \includegraphics[width=.5\textwidth]{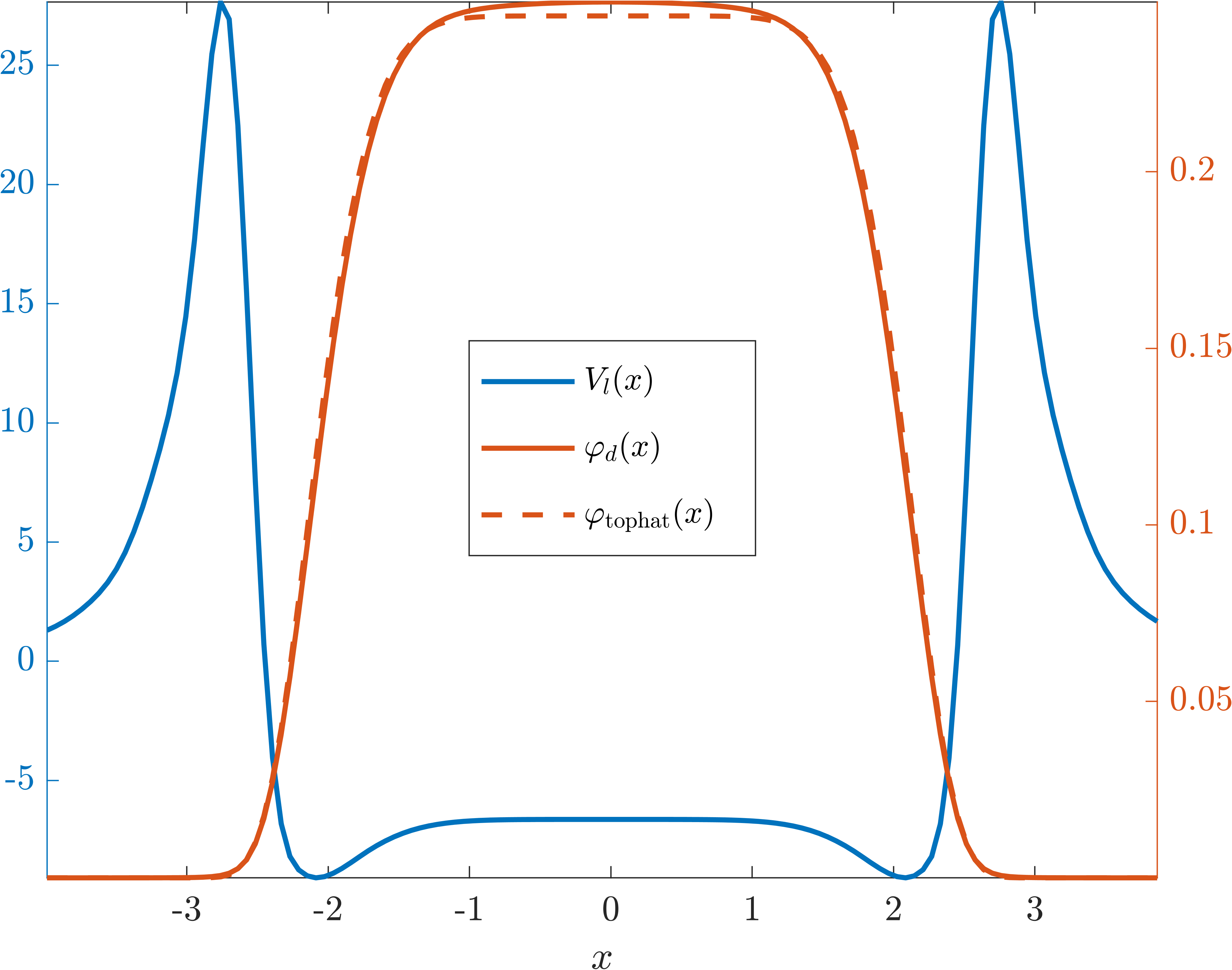}
    \caption{The top hat potential $V_l(x)$ which solves the inverse scattering problem~\eqref{eq:TopHatLsq} with top hat eigenfunction~\eqref{eq:TopHat}. The computed eigenfunction $\varphi_d(x)$ is in solid red.}
    \label{fig:TopHatPot}
\end{figure}

We now show how to solve two beam reshaping problems similiar to those originally considered by Kunkel and Leger~\cite{Kunkel}, with transverse dimension $n=$1, but by using the optimal control problem~\eqref{eq:BeamProb}. In the first problem, shown in Figure~\ref{fig:introsquare}, we transform the Pöschl-Teller eigenfunction
\begin{equation}
\varphi_0(x)=-\frac{1}{\sqrt{2}}\sech(x),
\end{equation}
which is the ground state of the potential 
\begin{equation}
\label{eq:PosclPot}
V_0(x)=-\frac{\sigma(\sigma+1)}{2}\sech^2(x)\end{equation}
when $\sigma=1$, into the "top hat" mode
\begin{equation}\label{eq:TopHat}
    \varphi_{\rm tophat}=Ae^{-ax^m},
\end{equation}
where $A$ is a normalization coefficient. For sake of computational demonstration, we choose $a=10^{-3}$ and $m=8$. The terminal potential $V_l(x)$ which has $\varphi_{\rm tophat}$ as its ground state mode is computed via the least squares problem
\begin{equation}\label{eq:TopHatLsq}
\min_{V_l(x)\in H_b^{1}\left(\mathbb{R}\right)}{J}=\min_{V_l(x)\in H_b^{1}\left(\mathbb{R}\right)}\frac{1}{2}\left\|\varphi_{\rm tophat}(x)-\varphi_d(x;V_l(x))\right\|_{L^2({\mathbb{R}})}^2
\end{equation}
subject to
\begin{equation}\label{eq:TopHatSchrod}
-\frac{1}{2}\partial_x^2\varphi_d(x)+V_l(x)\varphi(x)=\lambda_l\varphi_d(x).
\end{equation}
We show the resulting top hat potential $V_l(x)$ and eigenfunction $\varphi_d(x)$ from this procedure in Figure~\ref{fig:TopHatPot}. The computed eigenfunction $\varphi_d(x)$ is then used as a proxy for the true desired eigenfunction $\varphi_{\rm tophat}$ for the objective~\eqref{eq:BeamObj} of the optimal control problem.

With $V_l(x)$ computed, we address the corresponding beam reshaping problem. We reduce the search space of possible potentials by assuming they take the following form:
\begin{equation}\label{eq:assumedPot}
V(x,u(z),v(z))=u(z)V_0(x)+v(z)V_l(x),
\end{equation}
where $u(l)=v(0)=0$, and  $u(0)=v(l)=1$. This assumption on $V(x,z)$ slightly changes the optimality condition~\eqref{eq:controleq} such that the following equations
\begin{subequations}
\begin{alignat}{2}\label{eq:controlueq}
\gamma\partial_z^2{u}&=-\Re\left\langle p,V_0\psi\right\rangle_{L^2(\mathbb{R}^n)},\qquad &u(0)=u_0,\ u(l)=u_l,\\
\gamma\partial_z^2{v}&=-\Re\left\langle p,V_l\psi\right\rangle_{L^2(\mathbb{R}^n)},\qquad &v(0)=v_0,\ v(l)=v_l,
\end{alignat}
\end{subequations}
are now the appropriate Euler-Lagrange equations for the controls $u$ and $v$, while the state and costate equations~\eqref{eq:stateeq},~\eqref{eq:costateeq} remain unchanged.

\begin{figure}[htbp]
\begin{centering}
\subfigure[]{\includegraphics[width=.45\textwidth]{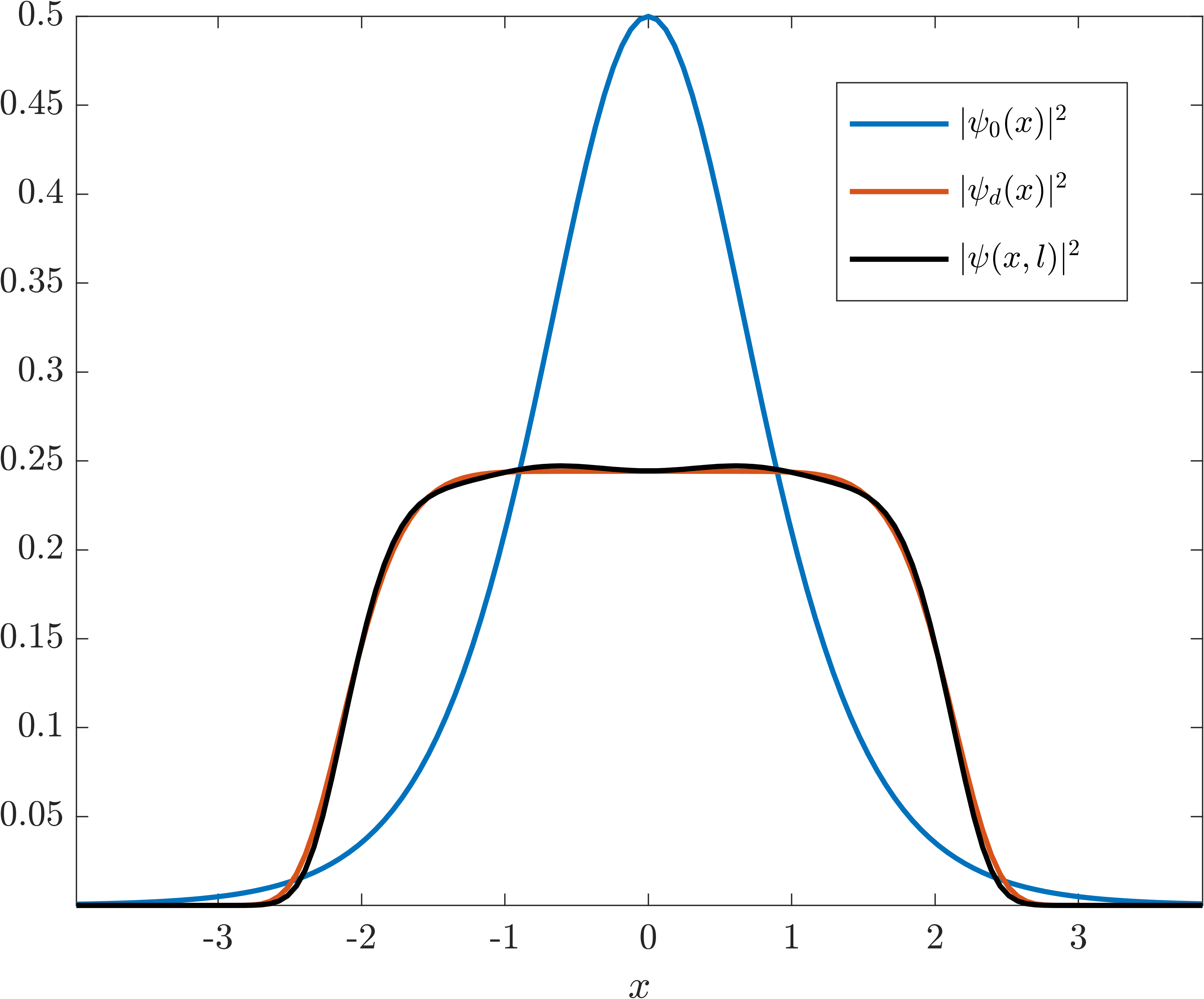}}
\subfigure[]{\includegraphics[width=.45\textwidth]{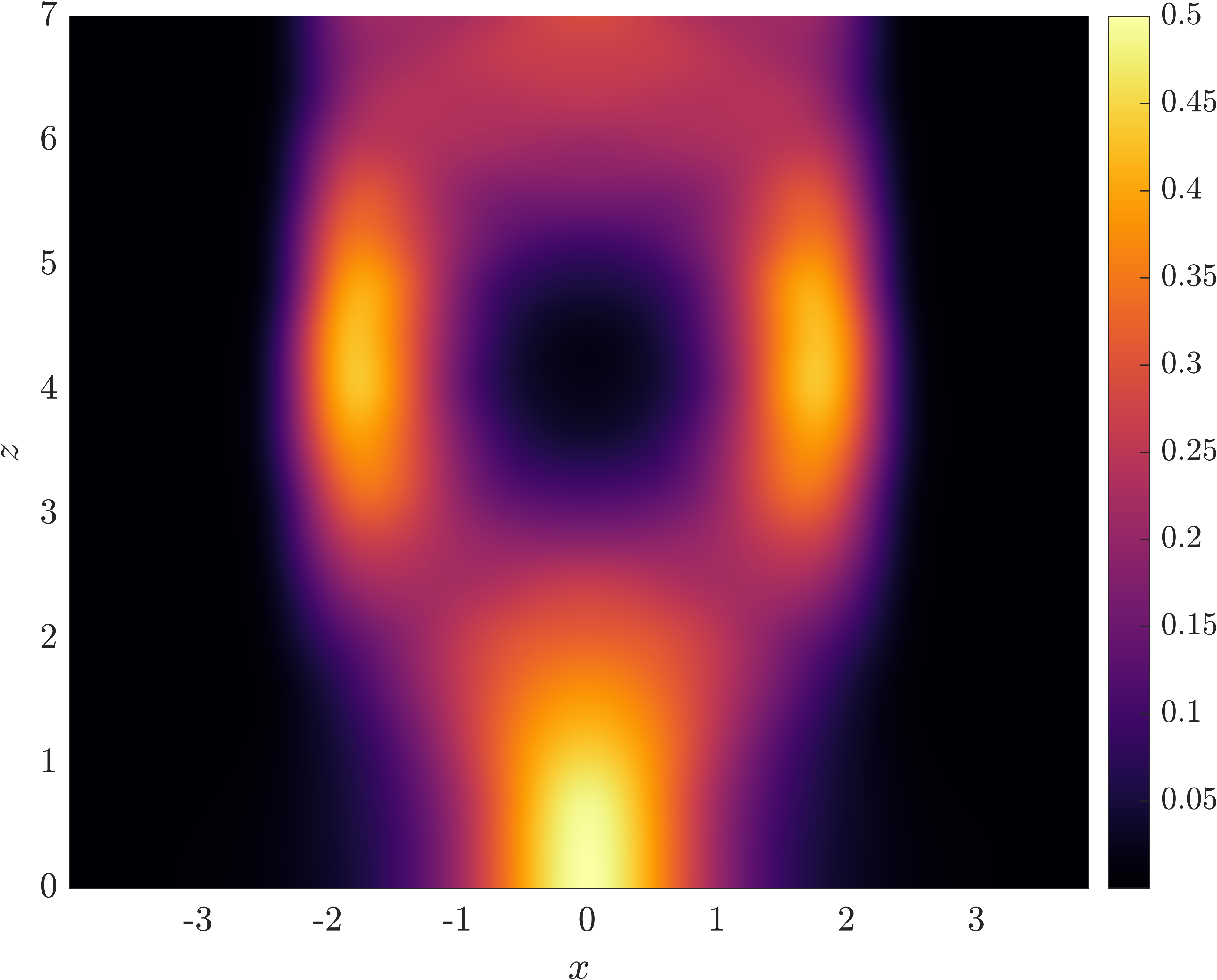}}
\subfigure[]{\includegraphics[width=.45\textwidth]{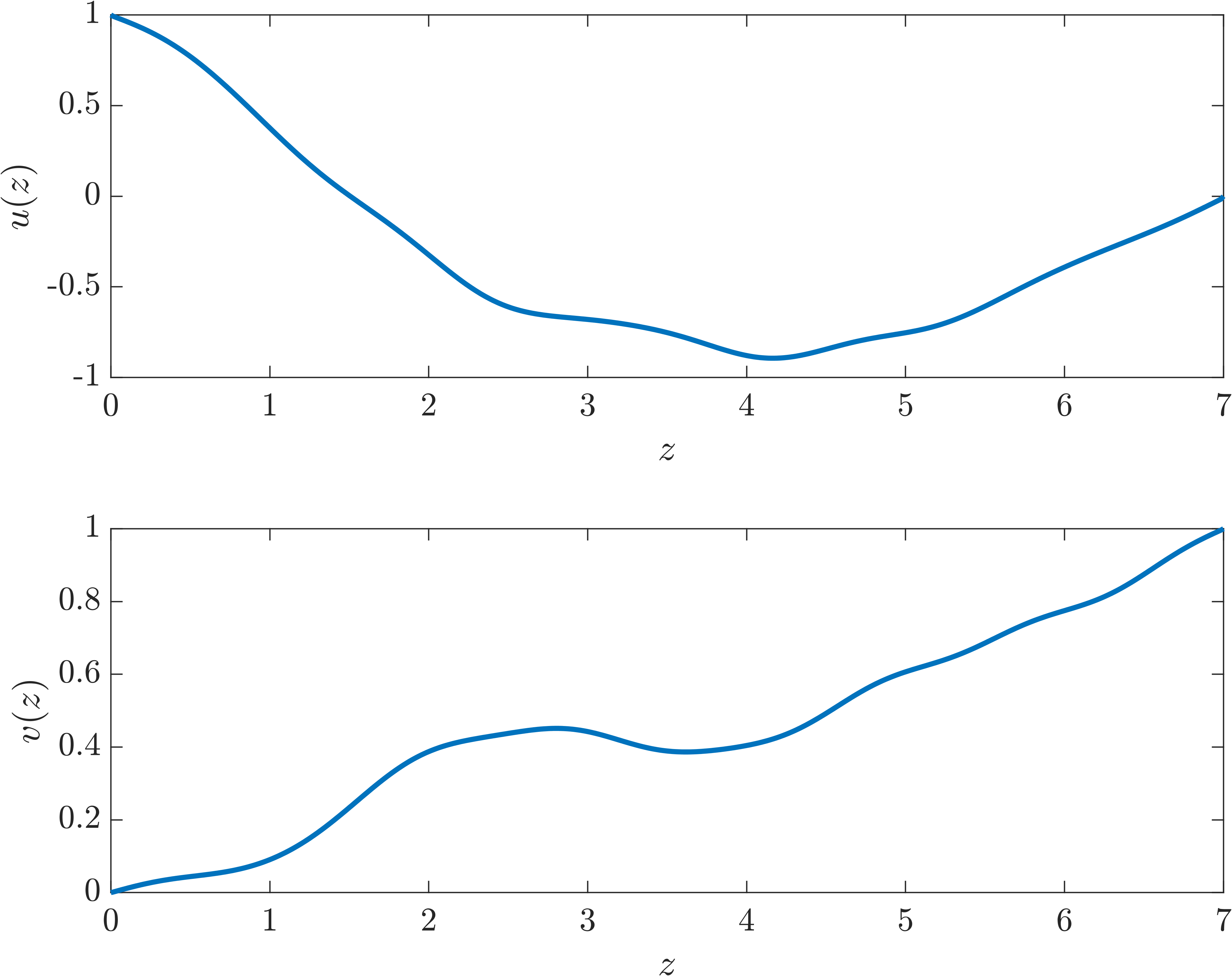}}
\subfigure[]{\includegraphics[width=.45\textwidth]{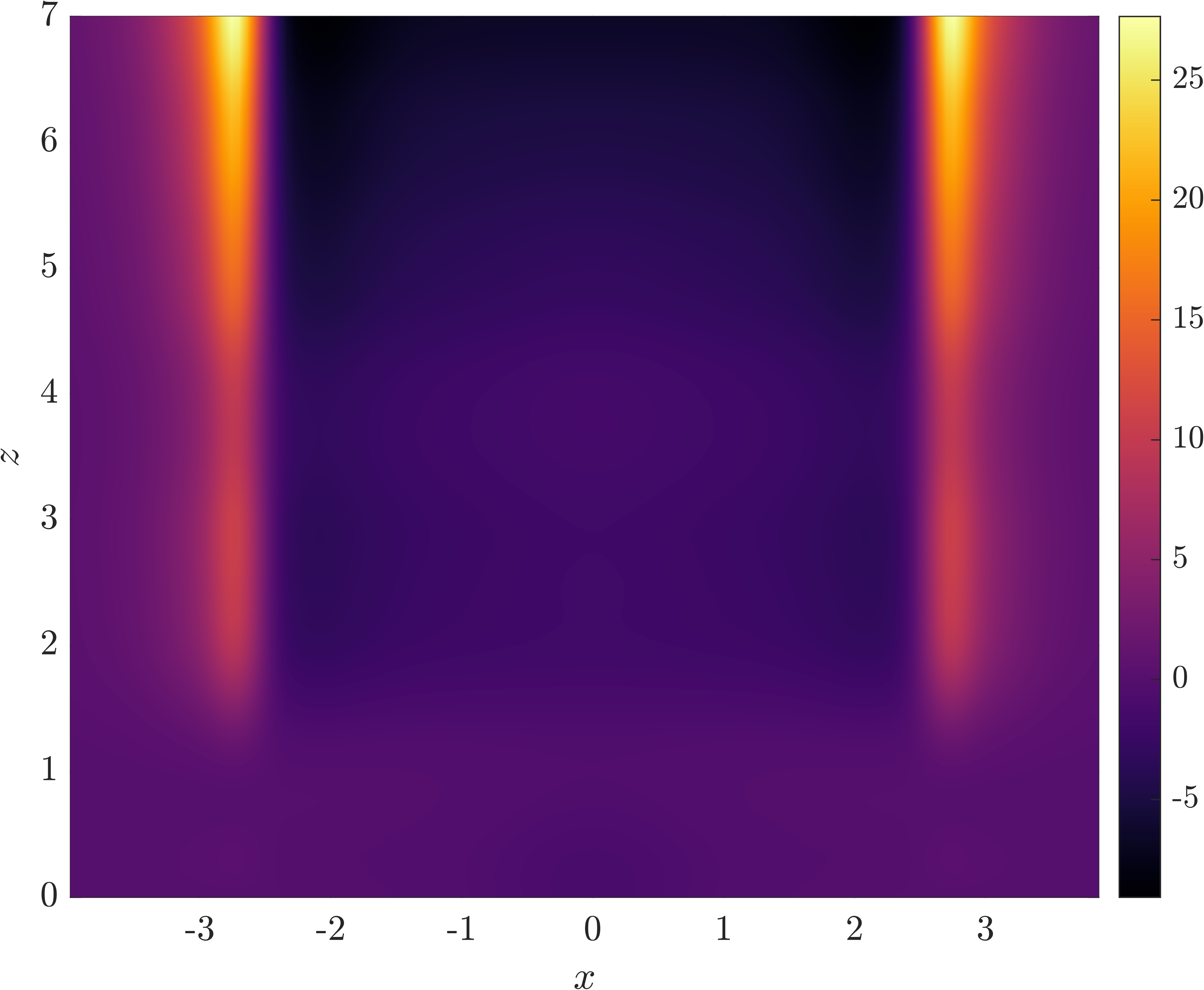}}
\caption{A numerical solution of the top hat problem. \textbf{(a)} The intensity profiles for the initial, desired, and final computed wavefunctions. \textbf{(b)} The axial evolution of the wavefunction intensity. \textbf{(c)} The computed controls $u(z)$ and $v(z)$ resulting from the hybrid method. \textbf{(d)} The optimal potential resulting from Panel (c) and the assumed form~\eqref{eq:assumedPot}.}
\label{fig:TopHatFig}
\end{centering}
\end{figure}

We show the results of the optimal control problem using the hybrid method of Section~\ref{section:beamnum} in Figure~\ref{fig:TopHatFig}. We set the Tikhonov parameter to $\gamma=10^{-6}$, fix $x\in[-5\pi,5\pi],\ z\in[0,7]$ and use a sine series together with a linear polynomial required by Equation~\eqref{eq:randclass}, i.e.,
\begin{subequations}\label{eq:PTCRAB}
\begin{align}
w_r(z)=\sum_{j=1}^{15}\frac{r_w}{j^2}\sin\left(\frac{j\pi z}{l}\right)+(w_l-w_0)\frac{z}{l}+w_0.
\end{align}
\end{subequations}
The amplitudes $r_w$ are random variables drawn uniformly from $[-1,1]$, and $w$ stands for either $u$ or $v$. We choose the coefficients $A_j=\frac{r_w}{j^2}$ to decay quadratically because the Fourier series of an absolutely continuous functions exhibits the same type of decay~\cite{Trefethen}. In this way, along with the relative smallness of the Tikhonov parameter, the search space for the optimal controls $u$ and $v$ is not severely restricted, yet candidate controls remain technologically feasible throughout each generation of DE and iteration of the projected gradient descent.

\subsection{The Beam Addition Problem}\label{section:Mprob}

Kunkel and Leger~\cite{Kunkel} consider the problem of merging several pulses into one, c.f. Figure~\ref{fig:introcombine}. To this end, we use an initial configuration of three seperated Pöschl-Teller potentials, each with $\sigma=1$, i.e.,  
\begin{align}\label{eq:beamcombinit}
V_0(x)&=-\left(\sech^2(x-a)+\sech^2(x+a)+\sech^2(x)\right),\\
\varphi_0(x)&=-\frac{1}{\sqrt{6}}\left(\sech(x-a)+\sech(x+a)+\sech(x)\right),
\end{align}
where the spacing parameter $a>0$. Although $\varphi_0(x)$ is not exactly an eigenfunction of $V_0(x)$, it approximates an eigenfunction with improving accuracy as $a$ is increased; we use $a=10$.

We emulate Kunkel and Leger's strategy of partitioning the the optimal control problem into two stages. In the context of this problem, we first perform an optimization on the interval $[0,30]$ where we use $V_0(x)$ as an initial potential and use the top hat potential of Figure~\ref{fig:TopHatPot} as the terminal potential. We then perform an optimization on the interval $[30,70]$ where the terminal data, i.e., the terminal potential and resulting terminal wavefunction, is used as initial data and the now terminal potential is given by a single Pöschl-Teller potential with $\sigma=3$. Both stages of the optimization are performed using the hybrid method on potentials of the form~\eqref{eq:assumedPot} with appropriate boundary conditions, in order to ensure continuity of potentials across $z=30$, and with the parameters $\gamma$ and $r_w$ the same as they were in Section~\ref{section:THprob}.

We further refine our results by relaxing the restriction of the search space from the assumed form~\eqref{eq:assumedPot} via a gradient descent on a wider space. That is, we perform a full two-dimensional gradient descent on the potential $V(x,z)$ resulting from the two-stage optimization. To compute the gradient in this case requires a solution of the Dirichlet problem
\begin{subequations}\label{eq:beamDirichlet}
\begin{align}
\label{eq:beampoisson}
    \nabla^2_{x,z}&\nabla_V\mathcal{J}=-\delta_V\mathcal{J},\\
\label{eq:beamhomo}
    &\nabla_V\mathcal{J}\big|_{\partial\Omega}=0,
    \end{align}
\end{subequations}
where the inhomogeneity is given by
\begin{equation}~\label{eq:sourceterm}
    -\delta_V\mathcal{J}=\gamma\nabla^2_{x,z}V+\Re\left\langle p,\psi\right\rangle_{L^2(\mathbb{R})},
\end{equation}
$\nabla^2_{x,z}$ is the Laplacian operator over $x$ and $z$, and $\partial\Omega$ is the boundary of the computational domain $[-15\pi,15\pi]\times[0,70]$. This is the GRAPE method, from Subsection~\ref{section:Local}, in the space $\dot{H}_0^1(\Omega)$. 

\begin{figure}[htbp]
\begin{centering}
\subfigure[]{\includegraphics[height=.4\textwidth]{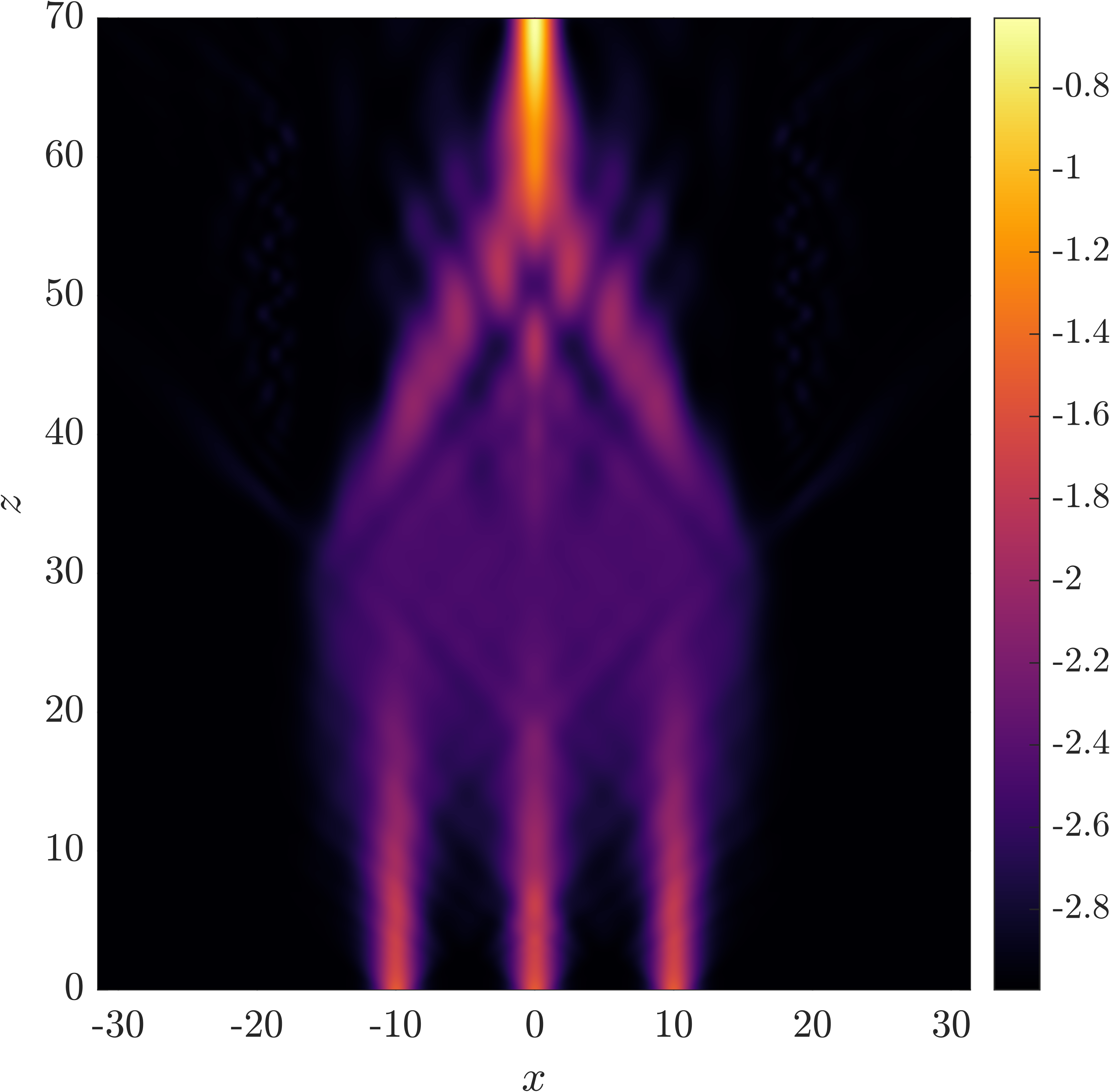}}
\subfigure[]{\includegraphics[height=.4\textwidth]{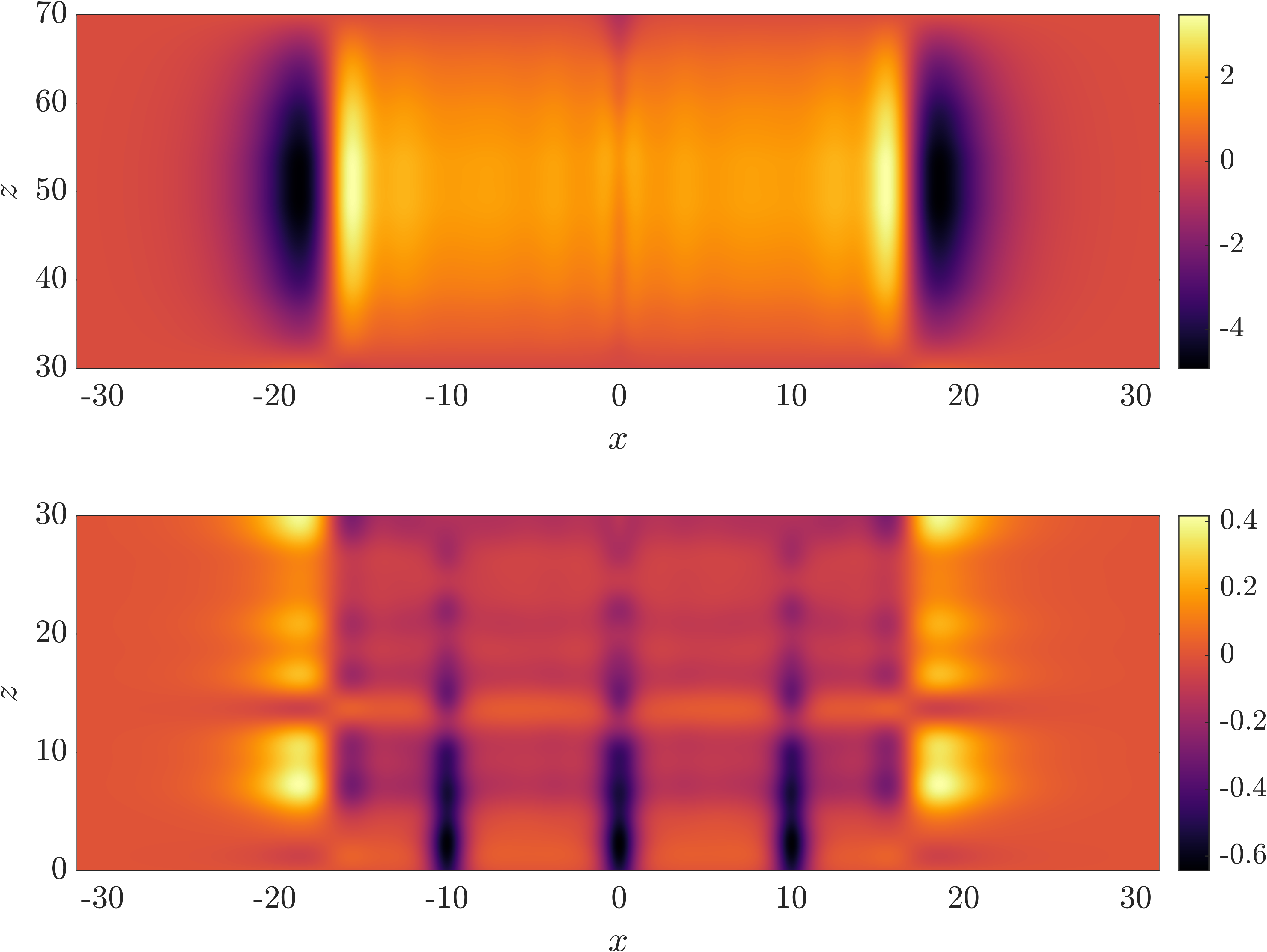}}
\caption{\textbf{(a)} The Schrödinger intensity distribution $|\psi|^2$, on a logarithmic scale. \textbf{(b)} The two stages of the computed optimal potential $V(x,z)$, shown separately because their ranges differ widely.}\label{fig:finalcombine1}
\end{centering}
\end{figure}

\begin{figure}[htbp]
\begin{centering}
\subfigure{\includegraphics[width=.45\textwidth]{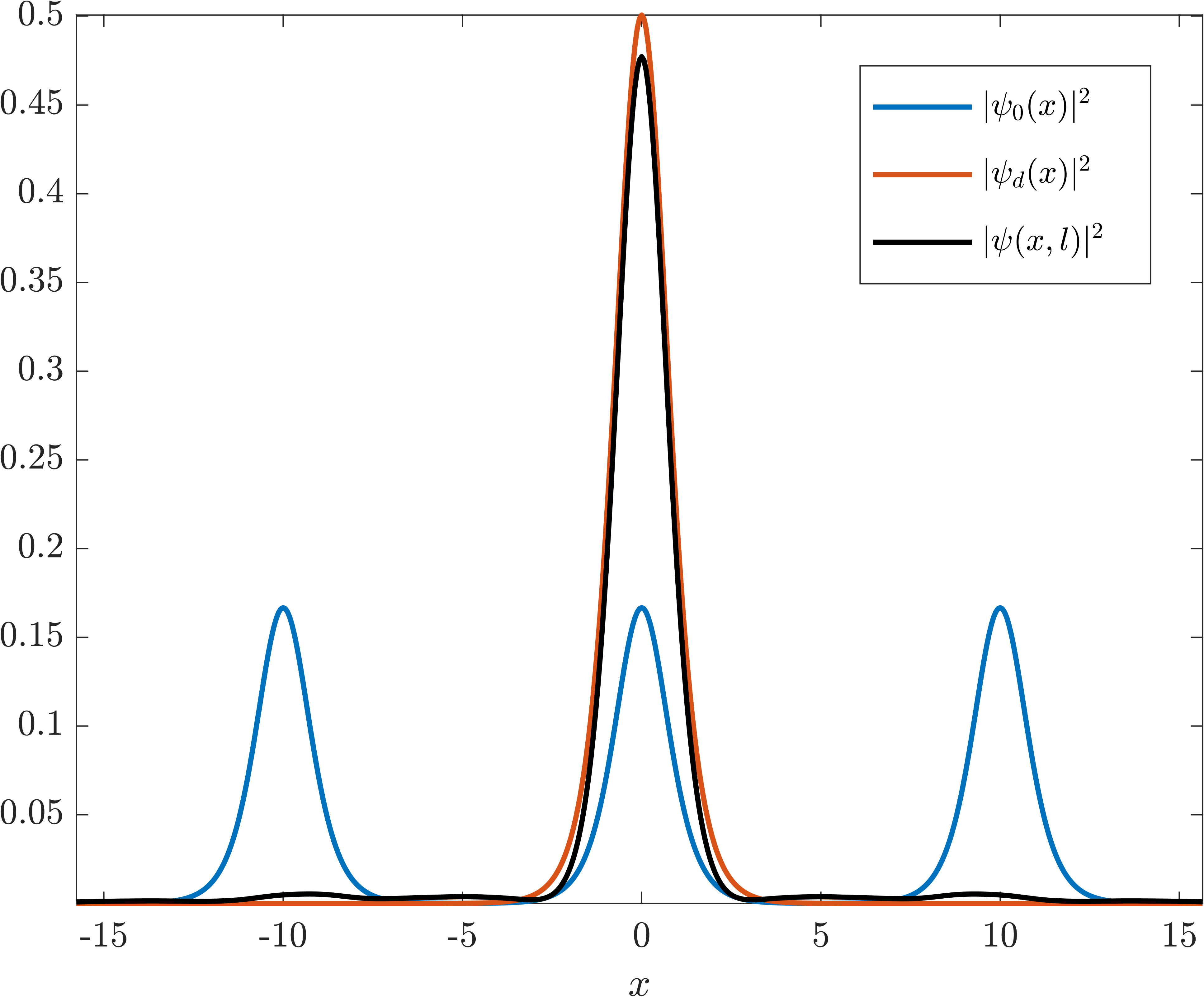}}
\caption{The initial, desired, and final computed intensity profiles corresponding to Figure~\ref{fig:finalcombine1}.}\label{fig:finalcombine2}
\end{centering}
\end{figure}

Note that the source term~\eqref{eq:sourceterm} in Poisson's equation~\eqref{eq:beampoisson} arises from Equation~\eqref{eq:controleq} and involves the computation of the Laplacian $\nabla^2_{x,z}$  which itself arises from the proper modification of the Tikhonov regularization in objective~\eqref{eq:BeamObj}, i.e., the cost now also runs over spatial dimension and penalizes large spatial derivatives of the reshaping potential $V$. We find this penalization is, on average, two orders of magnitude larger than penalizations which only run over the axial direction $z$, and so we decrease the Tikhonov parameter to $\gamma=10^{-8}$.  We show, in  Figures~\ref{fig:finalcombine1} and~\ref{fig:finalcombine2}, the final result of the GRAPE method in $\dot{H}_0^1(\Omega)$, after inputting the optimal controls computed through the two-stage hybrid optimization strategy.

\section{Conclusion and Future Work}
We have successfully applied optimal control theory  to design of GRIN fibers that reshape a beam of light into a desired shape. In order to thoroughly, yet efficiently, search the space of possible designs, we use a combination of a Galerkin reduction of the control, a projected gradient descent method, product separability of the reshaping potentials of the form $V(x,u(z),v(z))$, a partitioning of the control into stages, and finally gradient descents on a wider space for reshaping potentials of the full form $V(x,z)$. 

This methodology provides a systematic approach to the design process, but, of course, leaves further room for exploration. Moreover, the examples in this paper are proof of concepts in that we have only applied the methods to  waveguides with a single transverse dimension, whereas the phase retrieval method has now been applied to waveguides with two transverse direction~\cite{MinsterKunkel:20}. Future work may include extending the methods of this paper to higher dimensions. Fortunately, this extension is straightforward by virtue of the optimal control framework.

\section*{Acknowledgments}
We acknowledge the 2016-2017 program in optics at the University of Minnesota Institute for Mathematics and Applications (IMA) for connecting the authors with current research in optics, including that of the Leger group. We would also like to thank Alejandro Aceves and Braxton Osting for suggesting this work after attending that program. We also thank Richard Moore, John Federici, Louis Rizzo, David Shirokoff, James Leger, and Mint Kunkel for their helpful discussions, suggestions, and comments.

\bibliography{bibliography}

\end{document}